\documentclass[titlepage,11pt]{article}
\usepackage{latexsym}
\usepackage{amsfonts}
\usepackage{amsmath}
\newcommand\blackslug{\hbox{\hskip 1pt \vrule width 4pt height 8pt depth 1.5pt
        \hskip 1pt}}
\newcommand\bbox{\hfill \quad \blackslug \bigbreak}
\def\d{\hbox{-}}
\def\c{\hbox{-}\cdots\hbox{-}}
\def\l{,\ldots,}

\title{Induced subgraphs of graphs with large chromatic number.\\
III. Long holes}
\author{Maria Chudnovsky\thanks{Supported by NSF grant DMS-1265803.}\\
Princeton University, Princeton, NJ 08544, USA
\\
\\
Alex Scott\\
Mathematical Institute, University of Oxford, Oxford OX2 6GG, UK
\\
\\
Paul Seymour\thanks{Supported by ONR grant N00014-14-1-0084 and NSF
grant DMS-1265563.}\\
Princeton University, Princeton, NJ 08544}

\date{May 21, 2015; revised \today}

\newtheorem{thm}{}[section]

\newcommand{\Proof}{\noindent{\bf Proof.}\ \ }

\begin{document}
\maketitle
\begin{abstract}
We prove a 1985 conjecture of Gy\'arf\'as that for all $k,\ell$, every graph with sufficiently large chromatic number contains either 
a clique of cardinality more than $k$ or an induced cycle of length more than~$\ell$.
\end{abstract}

\section{Introduction}
All graphs in this paper are finite, and without loops or parallel edges. A {\em hole} in a graph $G$ is an induced subgraph
which is a cycle of length at least four, and an {\em odd hole} means a hole of odd length. (The {\em length} of a path
or cycle is the number of edges in it, and we sometimes call a hole of length $\ell$ an {\em $\ell$-hole}.)
In 1985, A. Gy\'arf\'as~\cite{gyarfas} made three beautiful and  well-known conjectures:

\begin{thm}\label{gyarfasconj0}
{\bf Conjecture:} For every integer $k>0$ there exists $n(k)$ such that every graph $G$ with no clique of cardinality
more than $k$ and no odd hole has chromatic number at most $n(k)$.
\end{thm}

\begin{thm}\label{gyarfasconj1}
{\bf Conjecture:} For all integers $k,\ell>0$ there exists $n(k,\ell)$ such that every graph $G$ with no clique of cardinality
more than $k$ and no hole of length more than $\ell$ has chromatic number at most $n(k,\ell)$.
\end{thm}

\begin{thm}\label{gyarfasconj}
{\bf Conjecture:} For all integers $k,\ell>0$ there exists $n(k,\ell)$ such that every graph $G$ with no clique of cardinality
more than $k$ and no odd hole of length more than $\ell$ has chromatic number at most $n(k,\ell)$.
\end{thm}

The third conjecture implies the first two, and remains open.
Virtually no progress was made on any of them until 2013, when two of us proved the first conjecture~\cite{oddholes}.
On the second and third, there was no progress at all, until we proved~\cite{pentagonal} that the second and third are 
both true when $k=2$ and $\ell =6$.
More recently two of us proved the third in full
when $k=2$~\cite{holeseq}. (In fact we proved much more; that for all $\ell\ge 0$, 
in every graph with large enough chromatic number
and no triangle, there is a sequence of holes of $\ell$ consecutive lengths). In this paper we prove the second; thus, our main result is:

\begin{thm}\label{mainthm}
For all integers $k,\ell>0$ there exists $c$ such that every graph $G$ with no clique of cardinality
more than $k$ and no hole of length more than $\ell$ has chromatic number at most $c$.
\end{thm}

We denote the chromatic number of a graph $G$ by $\chi(G)$.
If $X\subseteq V(G)$, the subgraph of $G$ induced on $X$ is denoted by $G[X]$,
and we often write $\chi(X)$ for $\chi(G[X])$. 

The proof of \ref{mainthm} is an extension of the method of~\cite{holeseq}. In particular, we proceed by 
induction on $k$, and so
we assume that \ref{mainthm} is true with $k$ replaced by $k-1$.
We will show in \ref{longholelemma} that if $G$ has no clique of cardinality
more than $k$ and no hole of length more than $\ell$, and has large chromatic number,
then there is a vertex $v_1$ such that the vertices within distance two of $v_1$ induce a subgraph with
(not quite so) large chromatic number. 
The set of neighbours of $v_1$ induces a subgraph with bounded chromatic number (because it contains
no clique of cardinality $k$), and so the vertices with distance exactly two from $v_1$ induce a subgraph, $G_2$ say,
with large chromatic number. By the same argument applied to $G_2$, there is a vertex $v_2$ of $G_2$ such that
the subgraph of $G_2$ induced on the vertices with distance exactly two from $v_2$
has large chromatic number. And so on, many times; we obtain a sequence of ``covers'', each covering the next, and
all covering some remaining subgraph which still has large chromatic number. The proof is by looking
closely at a long such sequence of covers.

\section{Multicovers}

Given a long sequence of covers, each covering the next as just explained, we can clean up the relationship between each
pair of them to make the relationship
between them as simple as possible, in a way that we explain later. It turns out that after the cleanup, there are 
two ways each pair of the 
covers might be related, and an application of Ramsey's theorem will give us a long subsequence where all the pairs are 
related the same way. Thus we need to extract something useful from a long sequence of covers pairwise related in the same 
way, where there are two possible cases for the ``way''. In this section we handle the first way; in that case we call 
the sequence of covers a ``multicover''.

If $X,Y$ are disjoint subsets of the vertex set of a graph $G$, we say 
\begin{itemize}
\item $X$ is {\em complete} to $Y$ if every vertex in $X$
is adjacent to every vertex in $Y$; 
\item $X$ is {\em anticomplete} to $Y$ if every vertex in $X$
nonadjacent to every vertex in $Y$; and
\item $X$ {\em covers} $Y$ if every vertex in $Y$ has a neighbour in $X$.
\end{itemize}
(If $X=\{v\}$ we say $v$ is complete to $Y$ instead of $\{v\}$, and so on.)
Let $x\in V(G)$, let $N$ be some set of neighbours of $x$, and let $C\subseteq V(G)$
be disjoint from $N\cup \{x\}$, such that $x$ is anticomplete to $C$ and $N$ covers $C$. In this situation
we call $(x,N)$ a {\em cover} of $C$ in $G$. For $C,X\subseteq V(G)$, a {\em multicover of $C$} in $G$
is a family $(N_x:x\in X)$
such that
\begin{itemize}
\item $X$ is stable;
\item for each $x\in X, (x,N_x)$ is a cover of $C$;
\item for all distinct $x,x'\in X$, $x'$ is anticomplete to $N_x$ (and in particular all the sets
$\{x\}\cup N_x$ are pairwise disjoint).
\end{itemize}
The multicover $(N_x:x\in X)$ is {\em stable} if each of the sets $N_x\;(x\in X)$ is stable. 
Let $(N_x:x\in X)$ be a multicover of $C$ in $G$. If $X'\subseteq X$, and $N_x'\subseteq N_x$ for each $x\in X'$, 
we say that $(N_x':x\in X')$ is {\em contained in} $(N_x:x\in X)$.

If $(N_x:x\in X)$ is a multicover of $C$, and $F$ is a subgraph of $G$ with $X\subseteq V(F)$ such that no vertex in 
$C\cup \bigcup_{x\in X}N_x$ belongs to
or has a neighbour in $V(F)\setminus X$, we say that $F$ is {\em tangent} to the multicover. 
We need to prove that if we are given a multicover $(N_x:x\in X)$ with $|X|$ large, of some set $C$ with $\chi(C)$ large, then there
is a multicover $(N_x':x\in X')$ of some $C'\subseteq C$, contained in $(N_x:x\in X)$, with $|X'|$ and $\chi(C')$ still large (but 
much smaller than before), and with
a certain desirable subgraph tangent, a ``tick''. 

Let $X\subseteq V(G)$ be stable. Let $a$ and $a_x\;(x\in X)$ be distinct members of $V(G)\setminus X$, such that
\begin{itemize}
\item $a$ is anticomplete to $X$;
\item $a_x$ is adjacent to $a,x$ and is anticomplete to $X\setminus \{x\}$, for each $x\in X$;
\end{itemize}
We call the subgraph of $G$ with vertex set $X\cup \{a\}\cup \{a_x:x\in X\}$ and edges $x\d a_x, a\d a_x$ for each $x\in X$
a {\em tick} on $X$ in $G$. This may not be an induced subgraph of $G$ because the vertices $a_x\;(x\in X)$ may be adjacent to one another in $G$.

For a graph $G$, we denote by $\omega(G)$ the cardinality of the largest clique of $G$, and if $X\subseteq V(G)$ we sometimes write $\omega(X)$ for
$\omega(G[X])$.
We need:

\begin{thm}\label{gettick}
For all $j,k,m,c,\kappa\ge 0$ there exist $m_j,c_j\ge 0$ with the following property.
Let $G$ be a graph with $\omega(G)\le k$, such that $\chi(H)\le \kappa$ for every induced subgraph $H$ of $G$ with $\omega(H)<k$.
Let $(N_x:x\in X)$ be a stable multicover in $G$ of some set $C$, such that $|X|\ge m_j$, $\chi(C)\ge c_j$,
and $\omega(\bigcup_{x\in X}N_x)\le j$.
Then there exist $X'\subseteq X$ with $|X'|\ge m$ and  $C'\subseteq C$ with $\chi(C')\ge c$ and a stable multicover $(N_x':x\in X')$
of $C'$ contained in $(N_x:x\in X)$, such that there is a tick in $G$ tangent to $(N_x':x\in X')$.
\end{thm}
\Proof We may assume that $k\ge 2$, for otherwise the result is vacuous. We proceed by induction on $j$, keeping $k,m,c,\kappa$ fixed. 
If $j=0$ then we may take $m_0 = c_0=1$ and the theorem holds vacuously; so we assume that $j>0$ and the result holds for $j-1$.  Thus $m_{j-1},c_{j-1}$ exist.
Let 
\begin{eqnarray*}
m_j &=& 2km m_{j-1}\\
d_2&=&m_j 2^{m_j}c_{j-1}+2^{m_j}c\\
d_1&=& d_2+m_j\kappa\\
d_0&=& k2^{m_j}d_1\\
c_j&=&d_0+k\kappa.
\end{eqnarray*}
We claim that $m_j,c_j$ satisfy the theorem.
Let  $G$, $(N_x:x\in X)$, and $C$ be as in the theorem, with
$|X|\ge m_j$ and $\chi(C)\ge c_j$,
such that $\omega(\bigcup_{x\in X}N_x)\le j$.
We may assume that $|X|=m_j$. 
Since $c_j>\kappa$, there is a clique $A\subseteq C$ with $|A|=k$. Let $C_0$ be the set of vertices in $C\setminus A$ with no 
neighbour in $A$; then since every vertex in $C\setminus C_0$ has a neighbour in $A$, and for each $a\in A$ its set of neighbours has
chromatic number at most $\kappa$ (because it includes no $k$-clique), it follows that $\chi(C\setminus C_0)\le k\kappa$, and so
$\chi(C_0)\ge c_j-k\kappa=d_0$.
\\
\\
(1) {\em There exist $a\in A$, and $X_1\subseteq X$ with $|X_1|\ge m_j/k$, and $C_1\subseteq C_0$ with 
$\chi(C_1)\ge d_1$, such that
for each $v\in C_1$ and each $x\in X_1$, there is a vertex in $N_x$ adjacent to $v$ and nonadjacent to $a$.}
\\
\\
For each $v\in C_0$ and each $x\in X$, $v$ has
a neighbour in $N_x$; and this neighbour is nonadjacent to some vertex in $A$, since $|A|=k=\omega(G)$. Thus there exists
$a_{v,x}\in A$ such that some vertex in $N_x$ is adjacent to $v$ and nonadjacent to $a_{v,x}$. There are only $k$ possible
values for $a_{v,x}$ as $x$ ranges over $X$, 
and so there exist $a_v\in A$ and $X_v\subseteq X$ with $|X_v|\ge |X|/k$, such that $a_{v,x}=a_v$ for all
$x\in X_v$. There are only $k$ possible
values for $a_v$; so there exist $a\in A$ and $C'\subseteq C_0$ with $\chi(C')\ge \chi(C_0)/k\ge 2^{m_j}d_1$, such that 
$a_v=a$ for all $v\in C'$. Thus for each $v\in C'$ there exists $X_v\subseteq X$ with $|X_v|\ge |X|/k$, such that $a_{v,x}=a$ for all
$x\in X_v$. There are at most $2^{m_j}$ possibilities for $X_v$; so there exists $C_1\subseteq C'$ with $\chi(C_1)\ge d_1$,
and $X_1\subseteq X$ with $|X_1|\ge m_j/k$, such that $X_v=X_1$ for all $v\in C_1$. This proves (1).

\bigskip

Let $a,X_1,C_1$ be as in (1). For each $v\in C_1$ and each $x\in X_1$, let $n_{x,v}\in N_x$ be adjacent to $v$ and nonadjacent to $a$.
For each $x\in X_1$ choose $a_x\in N_x$ adjacent to $a$.
Let $C_2$ be the set of all vertices in $C_1$ nonadjacent to each $a_x\;(x\in X_1)$; then $\chi(C_2)\ge \chi(C_1)-m_j\kappa\ge d_2$.
For each $y\in X_1$, let $C_y$ be the set of all $v\in C_2$ such that $n_{x,v}$ is adjacent to $a_y$,
for at least $m_{j-1}$ values of $x\in X_1\setminus \{y\}$.
Next, we show that we may assume that:
\\
\\
(2) {\em $\chi(C_y)\le c_{j-1}2^{m_j}$, for each $y\in X_1$.}
\\
\\
We will show that if (2) is false, then there is a multicover $(N'_x:x\in X')$ contained in 
$(N_x:x\in X)$ with $\omega(\bigcup_{x\in X'}N_x')\le j-1$, to which we can apply the 
the inductive hypothesis on $j$.
Suppose then that  $\chi(C_y)> c_{j-1}2^{m_j}$ for some $y\in X_1$.
For each $v\in C_y$, let $X_v\subseteq X_1\setminus \{y\}$ with $|X_v|=m_{j-1}$, such that $n_{x,v}$ is adjacent to $a_y$
for each $x\in X_v$. There are at most $2^{m_j}$ choices of $X_v$, and so there exist $C'\subseteq C_y$ and $X'\subseteq X_1\setminus \{y\}$
with $\chi(C')\ge \chi(C_y)2^{-m_j}\ge c_{j-1}$ and $|X'|=m_{j-1}$, such that $X_v=X'$ for all $v\in C'$. Let
$N'_x$ be the set of neighbours of $a_y$ in $N_x$, for each $x\in X'$; then $(N'_x:x\in X')$ is a multicover of $C'$.
Moreover, since every vertex in $\bigcup_{x\in X'}N'_x$ is adjacent to $a_y$, it follows that $\omega(\bigcup_{x\in X'}N'_x)<j$.
But then the result follows from the definition of $m_{j-1}, c_{j-1}$.
This proves (2).
\\
\\
(3) {\em There exist $C_3\subseteq C_2$ with $\chi(C_3)\ge c$ and 
$X_3\subseteq X_1$ with $|X_3|\ge m$, such that 
$n_{x,v}$ is nonadjacent to $a_y$ for all $v\in C_3$
and all distinct $x,y\in X_3$.}
\\
\\
Let $C'$ be the set of all $v\in C_2$ that are not in any of the sets $C_y\;(y\in X_1)$, that is, 
such that for each $y\in X_1$, there are fewer than $m_{j-1}$ values of $x\in X_1\setminus \{y\}$
such that $n_{x,v}$ is adjacent to $a_y$. From (2), it follows that 
$$\chi(C')\ge \chi(C_2)-m_j 2^{m_j}c_{j-1}\ge d_2- m_j 2^{m_j}c_{j-1}= 2^{m_j}c.$$
Let $v\in C'$; and let $G_v$ be the digraph with vertex set $X_1$ in which for distinct $x,y\in X_1$, $y$ is adjacent from $x$ in $G_v$
if $n_{x,v}$ is adjacent to $a_y$. It follows from the definition of $C_2$ that
every vertex of $G_v$ has indegree at most $m_{j-1}-1$. Consequently the undirected graph underlying $G_v$ has degeneracy at most 
$2m_{j-1}-2$, and therefore is $2m_{j-1}$-colourable. Thus there exists $X_v\subseteq X_1$ with $|X_v|\ge |X_1|/(2m_{j-1})$
such that no two members of $X_v$ are adjacent in $G_v$. There are at most $2^{m_j}$ choices of $X_v$, and so there exists $C_3\subseteq C'$
with $\chi(C_3)\ge \chi(C')2^{-m_j}\ge c$ and $X_3\subseteq X_1$ with 
$$|X_3|\ge |X_1|/(2m_{j-1})\ge m_j/(2km_{j-1})=m,$$ 
such that $X_v=X_3$ for all $v\in C_3$.
This proves (3).

\bigskip

For each $x\in X_3$, let $N_x'$ be the set of vertices in $N_x$ nonadjacent to each $a_y\;(y\in X_3)$. Thus $n_{x,v}\in N_x'$
for each $x\in X_3$ and $v\in C_3$. Hence $(N_x':x\in X_3)$ is a multicover of $C_3$ contained in $(N_x:x\in X)$. Moreover,
the subgraph consisting of $a$, the vertices $a_x\;(x\in X_3)$ and $X_3$, together with the edges $a\d a_x$ and $a_x\d x$ for each $x\in X_3$,
form a tick which is tangent to this multicover. This proves \ref{gettick}.~\bbox

\bigskip

Let us say an {\em impression} of $H$ in $G$ is a map $\eta$ with domain $V(H)\cup E(H)$, that maps $V(H)$ injectively into $V(G)$, and maps
each edge $e=uv$ of $H$ to a path of $G$ of length at least two joining the vertices $\eta(u),\eta(v)$; such that the set
$\{\eta(v):v\in V(H)\}$ is stable, and for every two edges $e,f$ of $H$ with no common end,
$V(\eta(e))$ is disjoint from and anticomplete to $V(\eta(f))$.
Its {\em order} is the maximum length of the paths $\eta(e)(e\in E(H))$.

By repeated application of \ref{gettick} with $j=|X|$, we can obtain many ticks on the same large subset $X'$ of $X$, disjoint except for $X'$ and with no edges
joining them disjoint from $X'$. (Note that vertices in the same tick with degree two in that tick may be adjacent in $G$, but otherwise
the subgraph formed by the union of the ticks is induced.) But such a ``tick cluster'' gives an impression of $K_{n,n}$ of order two, if
we take $n$ ticks clustered on a set $X'$ with $|X'|=n$. We deduce:

\begin{thm}\label{stableimpression}
Let $k,\kappa,n\ge 0$ be integers. Then there exist $m,c$ with the following property.
Let $G$ be a graph with $\omega(G)\le k$, such that there is no impression of $K_{n,n}$ in $G$ of order two, and such that
$\chi(H)\le \kappa$
for every induced subgraph $H$ of $G$ with $\omega(H)<k$. Then there is no stable multicover $(N_x:x\in X)$ in $G$
of a set $C$, such that $|X|\ge m$ and $\chi(C)\ge c$.
\end{thm}

Let us eliminate the ``stable'' hypothesis from \ref{stableimpression}.
\begin{thm}\label{impression}
Let $k,\kappa,n\ge 0$ be integers. Then there exist $m,c$ with the following property.
Let $G$ be a graph  with $\omega(G)\le k$, such that there is no impression of $K_{n,n}$ in $G$ of order two, and such that
$\chi(H)\le \kappa$
for every induced subgraph $H$ of $G$ with $\omega(H)<k$. Then there is no multicover $(N_x:x\in X)$ in $G$
of a set $C$, such that $|X|\ge m$ and $\chi(C)\ge c$.
\end{thm}
\Proof Let $m,c'$ satisfy \ref{stableimpression} (with $c$ replaced by $c'$).
Let $c=c'\kappa^{m}$. We claim that $m,c$ satisfy the theorem.
Let $G$ be as in the theorem, and suppose that $(N_x:x\in X)$ is a multicover in $G$
of a set $C$, such that $|X|\ge m$ and $\chi(C)\ge c$. We may assume that $|X|=m$.
For each $x\in X$, the subgraph induced on $N_x$ is $\kappa$-colourable; choose some such colouring, with colours $1\l \kappa$, for each $x$.
For each $v\in C$, let $f_v:X\rightarrow \{1\l \kappa\}$ such that for each $x\in X$, some neighbour of $v$ in $N_x$
has colour $f_v(x)$. There are only $\kappa^{|X|}$ possibilities for $f_v$, so there is a function $f:X\rightarrow \{1\l \kappa\}$
and a subset $C'\subseteq C$ with $\chi(C')\ge \chi(C)\kappa^{-|X|}\ge c'$, such that $f_v=f$ for all $v\in C'$. For each $x\in X$,
let $N'_x$ be the set of vertices in $N_x$ with colour $f(x)$; then $(N'_x:x\in X)$ is a stable multicover of $C'$, and the result follows
from the choice of $m,c'$. This proves \ref{impression}.~\bbox

If $G$ admits an impression of $K_{n,n}$, then $G$ has a hole of length at least $2n$. We deduce
\begin{thm}\label{nomulticover}
Let $k,\kappa,\ell\ge 0$ be integers. Then there exist $m,c$ with the following property.
Let $G$ be a graph with no hole of length at least $\ell$, such that $\omega(G)\le k$, and $\chi(H)\le \kappa$
for every induced subgraph $H$ of $G$ with $\omega(H)<k$. Then there is no multicover $(N_x:x\in X)$ in $G$
of a set $C$, such that $|X|\ge m$ and $\chi(C)\ge c$.
\end{thm}

We remark that with a little more work, we can prove a version of \ref{gettick}, and of \ref{nomulticover} below,
 which just assumes there is no odd hole of 
length at least $\ell$, instead of assuming there is no hole of length at least $\ell$.
The proof is, roughly: use the argument above to get a large tick cluster, all tangent to a multicover $(N_x:x\in X)$
of some set $C$, with $|X|$ and $\chi(C)$ large. Use Ramsey's theorem repeatedly, to arrange that for each tick,
its ``knees'' are stable (shrinking $X$ to some smaller set); and then choose an odd path between
two vertices $x,x'\in X$ via a vertex in $N_x$, a vertex in $N_{x'}$, and an $\omega(G)$-clique in $C$.
We omit the details.

\section{Cables}

Now we return to the long sequence of covers mentioned at the start of the previous section. The goal of this section
is just to
introduce some terminology, describing precisely what results after the clean-up process (but before the application of 
Ramsey's theorem), and then to carry out the application of Ramsey's theorem. 

Let $X\subseteq V(G)$ be a clique. If $|X|=k$ we call $X$ a {\em $k$-clique.}
We denote by $N^1_G(X)$ the set of all vertices in $V(G)\setminus X$ that are complete to $X$;
and by $N^2_G(X)$ the set of all vertices in $V(G)\setminus X$ with a neighbour in $N^1(X)$ and with no neighbour in $X$.
When $X=\{v\}$ we write $N^i_G(v)$ for $N^i_G(X)$ ($i = 1,2$).

Let $G$ be a graph and let $t\ge 0$ and $h\ge 1$ be integers. An {\em $h$-cable} of {\em length $t$} in $G$ consists of:
\begin{itemize}
\item $t$ $h$-cliques $X_1\l X_t$, pairwise disjoint and anticomplete;
\item for $1\le i\le t$, a subset $N_i$ of $N^1_G(X_i)$, such that the sets $N_1\l N_t$ are pairwise disjoint;
\item for $1\le i\le t$, disjoint subsets $Y_{i,t}$ and $Z_{i,i+1}\l Z_{i,t}$ of $N_i$; and
\item a subset $C\subseteq V(G)$ disjoint from $X_1\cup\cdots\cup X_t\cup N_1\cup\cdots\cup N_t$
\end{itemize}
satisfying the following conditions.
\begin{enumerate}
\item[{\bf (C1)}] For $1\le i\le t$, $Y_{i,t}$ covers $C$, and $C$ is anticomplete to $Z_{i,j}$ for $i+1\le j\le t$, and 
$C$ is anticomplete to $X_i$.
\item[{\bf (C2)}] For $i<j\le t$, $X_i$ is anticomplete to $N_j$.
\item[{\bf (C3)}] For all $i<j\le t$, every vertex in $Z_{i,j}$ has a non-neighbour in $X_j$.
\item[{\bf (C4)}] For $i<j<k\le t$, $Z_{i,j}$ is anticomplete to $X_k\cup N_k$.
\item[{\bf (C5)}] For all $i<j\le t$, either
\begin{itemize}
\item some vertex in $X_j$ is anticomplete to $Y_{i,t}$, and $Z_{i,j}=\emptyset$, or
\item $X_j$ is complete to $Y_{i,t}$, and $Z_{i,j}$ covers $N_j$.
\end{itemize}
\end{enumerate}

We call $C$ the {\em base} of the $h$-cable, and say $\chi(C)$ is the {\em chromatic number} of the $h$-cable.
Given an $h$-cable in this notation, let $I\subseteq \{1\l t\}$;
then the cliques $X_i\;(i\in I)$, the sets $N_{i}\; (i\in I)$, the sets $Z_{i,j}\; (i,j\in I)$, the sets $Y_{i,t}\; (i\in I)$
and $C$ (after appropriate renumbering) define an $h$-cable of length $|I|$. We call this a {\em subcable}.

Thus there are two types of pair $(i,j)$ with $i<j\le t$, and later we will apply Ramsey's theorem on these pairs to get a large subcable
where all the pairs have the same type.
Consequently, two special kinds of $h$-cables are of interest ($t$ is the length in both cases): 
\begin{itemize}
\item $h$-cables of {\em type 1},
 where for all $i<j\le t$, some vertex in $X_j$ has no neighbours in $Y_{i,t}$, and $Z_{i,j}=\emptyset$; and
\item $h$-cables of {\em type 2}, where for all $i<j\le t$,
$X_j$ is complete to $Y_{i,t}$, and $Z_{i,j}$ covers $N_j$.
\end{itemize}

From \ref{nomulticover} we deduce:

\begin{thm}\label{usetype1}
For all $k,\kappa,n\ge 0$ and $h\ge 1$, there exist $t,c\ge 0$ with the following property.
Let $G$ be a graph with  $\omega(G)\le k$, such that there is no impression of $K_{n,n}$ in $G$ of order two, and 
$\chi(H)\le \kappa$
for every induced subgraph $H$ of $G$ with $\omega(H)<k$. 
Then $G$ admits no $h$-cable of type 1 and length $t$ with chromatic number more than $c$.
\end{thm}
\Proof Choose $m,c$ to satisfy \ref{impression}. 
By Ramsey's theorem there exists $t$ such that for every partition of the edges of $K_t$ into
$h$ sets, there is an $m$-clique of $K_t$ for which all edges joining its vertices are in the same set.
We claim that $t,c$ satisfy the theorem. 

For let $G$ be as in the theorem, and
suppose that $G$ admits an $h$-cable of type 1 and length $t$ with chromatic number more than $c$.
In the usual notation for $h$-cables, fix an ordering of the members of $X_i$ for each $i$; thus we may speak of the 
$r$th member of $X_i$ for $1\le r\le h$. For each pair $(i,j)$ with $i<j\le t$, let $f(i,j) = r$ where the $r$th member
of $X_j$ has no neighbours in $Y_{i,t}$. From the choice of $t$, there exist $I\subseteq \{1\l t\}$ with $|I|=m$ 
and $r\in \{1\l h\}$
such that $f(i,j)=r$ for all $i,j\in I$ with $i<j$. For each $j\in I$, let $x_j$ be the $r$th member of $X_j$.
Then the sets $(x_j,N_j)\;(j\in I)$ form a multicover of $C$, which is impossible by \ref{impression}.
This proves \ref{usetype1}.~\bbox

\bigskip

We need an analogue for cables of type 2, but it needs an extra hypothesis.
On the other hand, we only need to assume that there is no hole of length exactly $\ell$.
\begin{thm}\label{usetype2}
Let $\tau\ge 0$, $\ell\ge 5$ and $h\ge 1$, and let $G$ be a graph with no $\ell$-hole, such that
$\chi(N^2(X))\le \tau$ for every $(h+1)$-clique $X$ of $G$.
Then $G$ admits no $h$-cable of type 2 and length $\ell-3$ with chromatic number more than $(\ell-3)\tau$.
\end{thm}
\Proof
Let $t=\ell-3$, let $G$ be as in the theorem, and
suppose that $G$ admits an $h$-cable of type 2 and length $t$ with chromatic number more than $t\tau$.
In the usual notation, choose $z_t\in Y_{t,t}$ (this is possible by {\bf (C1)}), and choose 
$z_{t-1}\in Z_{t-1,t}$ adjacent to $z_t$ (this is possible since the cable has type 2).
Since $z_{t-1}\in Z_{t-1,t}$, it has a non-neighbour $x_t\in X_t$, by {\bf (C3)}. 
Neither of $x_t,z_t $ has a neighbour in
$Z_{i,i+1}$ for $1\le i \le t-2$, by {\bf (C4)}. 
Now $z_{t-1}$ has a neighbour $z_{t-2}\in Z_{t-2,t-1}$, since the cable has type 1; and similarly 
for $i=t-3\l 1$
let $z_i\in Z_{i,i+1}$ be a neighbour of $z_{i+1}$. It follows that 
$$z_1\d z_2\c z_{t-1}\d z_t\d x_t$$
is an induced path.

For $1\le i\le t$, let $C_i$ be the set of vertices $v\in C$ such that some vertex in $Y_{1,t}$ is adjacent to both $v,z_{i}$.
Since $X_i$ is complete to $Y_{1,t}$ (since the cable has type 2), 
it follows that $C_i\subseteq N^2_G(X_i\cup \{z_i\})$; and since
$X_i\cup \{z_i\}$ is an $(h+1)$-clique, it follows from the hypothesis that $\chi(C_i)\le \tau$. Thus the union
$C_1\cup\cdots \cup C_{t}$ has chromatic number at most $t\tau$; and since $\chi(C)>t\tau$, 
there exists $u\in C$ not
in any of the sets $C_i\:(1\le i\le t)$. Choose $v\in Y_{1,t}$ adjacent to $u$ (this is possible by {\bf (C1)});
then $v$ is not adjacent to
any of $z_1\l z_{t}$, by definition of $C_1\l C_{t}$. Choose $x_1\in X_1$; then 
$$v\d x_1\d z_1\d z_2\c z_{t-1}\d z_t\d x_t\d v$$
is a hole of length $t+3=\ell$, a contradiction. This proves \ref{usetype2}.~\bbox

\bigskip

From \ref{usetype1}, \ref{usetype2} and Ramsey's theorem, we deduce:

\begin{thm}\label{usecable}
For all $k,\kappa,\tau,\ell\ge 0$ and $h\ge 1$, there exist $t,c\ge 0$ with the following property.
Let $G$ be a graph such that:
\begin{itemize}
\item $G$ has no hole of length at least $\ell$; 
\item $\omega(G)\le k$; 
\item $\chi(H)\le \kappa$ for every induced subgraph $H$ of $G$ with $\omega(H)<k$; and
\item $\chi(N^2(X))\le \tau$ for every $(h+1)$-clique $X$ of $G$.
\end{itemize}
Then every $h$-cable in $G$ of length $t$ has chromatic number at most $c$.
\end{thm}

\section{Clique control}

In this section we explain the clean-up process that we 
plan to apply to the long sequence of covers; but before 
that, we need another concept, ``clique control''.

Let $\mathbb{N}$ denote the set of nonnegative integers, let $\phi:\mathbb{N}\rightarrow \mathbb{N}$ be a nondecreasing function, and
let $h\ge 1$ be an integer. We say a graph $G$ is {\em $(h,\phi)$-clique-controlled} if for every induced subgraph $H$ of $G$ and
every integer $n\ge 0$, if $\chi(H)>\phi(n)$ then there is an $h$-clique $X$ of $H$ such that $\chi(N^2_H(X))>n$.
Intuitively, this means that in every induced subgraph $H$ of large chromatic number, there is an $h$-clique
$X$ with $N^2_H(X)$ of large chromatic number; the function $\phi$ is just a way of making ``large'' precise.

The following contains the clean-up process, somewhat disguised (rather than choose the whole sequence of covers
and then clean up all the pairs of its terms separately, it is more convenient to grow the sequence term by term
cleaning up all pairs involving the new term at each step).

\begin{thm}\label{getcable}
Let $t,c,\tau,\kappa\ge 0$ and $h>0$, and let $\phi:\mathbb{N}\rightarrow\mathbb{N}$ be nondecreasing.
Then there exists $c'$ with the following property.
Let $G$ be a graph such that
\begin{itemize}
\item $\chi(N^1(v))\le \kappa$ for every $v\in V(G)$;
\item $G$ is $(h,\phi)$-clique-controlled; and
\item $\chi(N^2(X))\le \tau$ for every $(h+1)$-clique $X$ of $G$.
\end{itemize}
If $\chi(G)> c'$ then
$G$ admits an $h$-cable of length $t$ with chromatic number more than $c$.
\end{thm}
\Proof Let $\sigma_t=\max(c,\tau+h\kappa)$, and for $s=t-1\l 0$ let
$$\sigma_s=\max(2^s\phi((h+1)^{s}\sigma_{s+1}),\tau+h\kappa).$$
Let $c'= \sigma_0$.  We claim that $c'$ satisfies the theorem.

Let $G$ be a graph satisfying the hypotheses of the theorem, and therefore with $\chi(G)> c'$.
Consequently $G$ admits an $h$-cable of length $0$ with chromatic number more than $\sigma_0$. We claim that
for $s=1\l t$, $G$ admits an $h$-cable of length $s$ with chromatic number more than $\sigma_s$. For suppose the
result holds for some $s<t$; we will prove it also holds for $s+1$.

Thus, $G$ admits an $h$-cable of length $s$ with chromatic number more
than $\sigma_s$.
In the usual notation, let $C$ be the base of the $h$-cable. For each $v\in C$ and $1\le i\le s$, let $C_{i,v}$ be the set of vertices $u\in C\setminus \{v\}$
nonadjacent to $v$, such that some vertex in $Y_{i,s}$
is adjacent to both $u,v$. Let $f_{i,v}=1$ if $\chi(C_{i,v})>\tau+h\kappa$, and $f_{i,v}=0$ otherwise. There are only $2^s$ possibilities
for the sequence $f_{1,v}\l f_{s,v}$, so there is a subset $C_1\subseteq C$ with $\chi(C_1)\ge  2^{-s}\chi(C)>2^{-s}\sigma_s$
and a $0,1$-sequence $f_1\l f_s$
such that $f_{i,v}=f_i$ for $1\le i\le s$ and all $v\in C_1$.
For $0\le i\le s$ let $d_i=(h+1)^{s-i}\sigma_{s+1}$.
Let $H=G[C_1]$; then
since $2^{-s}\sigma_s\ge \phi(d_0)$, there is an $h$-clique $X_{s+1}$ of $H$ such that $\chi(D_0)> d_0$, where $D_0=N^2_H(X_{s+1})$.
Let $N_{s+1} = Y_{s+1,s+1}=N^1_H(X_{s+1})$.

For $1\le i\le s$, we define $Y_{i,s+1},Z_{i,s+1}\subseteq Y_{i,s}$ and $D_i\subseteq D_{i-1}$ as follows.
Assume that we have defined $D_{i-1}$, and $\chi(D_{i-1})>d_{i-1}$.
Let $W$ be the set of vertices in $Y_{i,s}$ that are complete to $X_{s+1}$, and for each $x\in X_{s+1}$, let
$U_x$ be the set of vertices in $D_{i-1}$
with a neighbour in $Y_{i,s}$ that is nonadjacent to $x$. If $\chi(U_x)> d_i$ for some $x\in X_{s+1}$, let $D_i=U_x$,
let $Y_{i,s+1}$ be the set of all vertices in $Y_{i,s}$ that are nonadjacent to $x$, and let $Z_{i,s+1}=\emptyset$.
Let us call this ``case 1''.

Thus we assume that $\chi(U_x)\le d_i$ for each $x\in X_{s+1}$; and so $\bigcup_{x\in X_{s+1}}U_x$ has chromatic number
at most $hd_i$. Let $D_i=D_{i-1}\setminus \bigcup_{x\in X_{s+1}}U_x$; then $\chi(D_i)> d_{i-1}-hd_i=d_i$.
For each vertex in $D_i$, all its neighbours in $Y_{i,s}$ belong to $W$. In particular, let $x\in X_{s+1}$; then $C_{i,x}$ (defined earlier)
has chromatic number more than
$$d_i\ge \sigma_{s+1}\ge \tau+h\kappa,$$
and so $f_{i,x}=1$. Since $x\in C_1$, it follows that
$f_i=1$, and so $\chi(C_{i,v})>\tau+h\kappa$ for each $v\in C_1$.

Now let $v\in N_{s+1}$. If $u\in C$, and $u$ has no neighbour in $X_{s+1}\cup \{v\}$, and some vertex in $W$ is adjacent to
both $u,v$, then $u\in N^2_G(X_{s+1}\cup\{v\})$; and so the set of all such $u$ has chromatic number at most $\tau$.
On the other hand, the set of $u\in C$ with a neighbour in $X_{s+1}$ has chromatic
number at most $h\kappa$, since for each $x\in X_{s+1}$
its set of neighbours has chromatic number at most $\kappa$. Consequently the set of vertices in $C$
that are nonadjacent to $v$ and adjacent to a neighbour of $v$ in $W$ has chromatic number at most $\tau+h\kappa$.
Since $\chi(C_{i,v})>\tau+h\kappa$, it follows that there exists $u\in C_{i,v}$ such that no neighbour of $v$ in $W$
is adjacent to $u$. From the definition of $C_{i,v}$, it follows that $v$ has a neighbour in $Y_{i,s}\setminus W$.

Since this is true for every vertex $v\in N_{s+1}$,
we may define $Y_{i,s+1}=W$ and $Z_{i,s+1}=Y_{i,s}\setminus W$, and it follows that $Z_{i,s+1}$ covers $N_{s+1}$. 
This completes the definition of
$Y_{i,s+1},Z_{i,s+1}$ and $D_i$. Let us call this ``case 2''.

In either case, $\chi(D_s)>d_s$, and 
we claim that $X_1\l X_{s+1}$, the sets $N_1\l N_{s+1}$, the sets $Z_{i,j}$ for $1\le i<j\le s+1$,
the sets $Y_{i,s+1}$ for $1\le i\le s+1$, and $D_s$, define
an $h$-cable of length $s+1$ and chromatic number more than $d_s$.
To see this, we  must verify {\bf(C1)}--{\bf(C5)}. 

For {\bf(C1)}, since $D_s$ is anticomplete to $Z_{i,j}$ for $i+1\le j\le s$, and $D_s$ is anticomplete to $X_i$ for $1\le i\le s$,
and $D_s$ is anticomplete to $X_{s+1}$ from its definition,
it is enough to show that for $1\le i\le s+1$,
$Y_{i,s+1}$ covers $D_{s}$, and if $i\le s$ then $D_s$ is anticomplete to $Z_{i,s+1}$.
Suppose first that $i=s+1$. Since $D_s\subseteq D_0=N^2_H(X_{s+1})$ and
$Y_{s+1,s+1}=N^1_H(X_{s+1})$, it follows that $Y_{s+1,s+1}$ covers $D_{s}$, 
so the first claim holds; and the second holds vacuously. Thus we may assume
that $1\le i\le s$. Assume that case 1 applies, and let $x$ be as in case 1. We recall that $D_i=U_x$,
$Y_{i,s+1}$ is the set of all vertices in $Y_{i,s}$ that are nonadjacent to $x$, and $Z_{i,s+1}=\emptyset$. Consequently
$Y_{i,s+1}$ covers $U_x$ (from the definition if $U_x$) and hence covers $D_s\subseteq D_i=U_x$, from the definition of $U_x$, so the first claim holds.
Now $D_s$ is anticomplete to $Z_{i,s+1}$ since the latter 
is empty; so the second claim holds. 
This proves {\bf(C1)} in case 1. Now we assume that case 2 applies. 
With notation as in case 2, we recall that  $D_i=D_{i-1}\setminus \bigcup_{x\in X_{s+1}}U_x$, $Y_{i,s+1}=W$ and 
$Z_{i,s+1}=Y_{i,s}\setminus W$. For every vertex in $D_i$, all its neighbours in $Y_{i,s}$ belong to $W$, and so $Z_{i,s+1}$ is anticomplete to $D_s\subseteq
D_i$, and the second claim holds; and since every vertex in $D_i$ has
such a neighbour, it follows that $Y_{i,s+1}$ covers $D_{s}\subseteq D_i$, so the first claim holds.
This completes the proof of {\bf(C1)}.

For {\bf(C2)}, it suffices to show that for $i< s+1$, $X_i$ is anticomplete to $N_{s+1}$. But this is true since $N_{s+1}\subseteq C$
and $X_i$ is anticomplete to $C$. This proves {\bf(C2)}.

For {\bf(C3)}, it suffices to show that for all $i< s+1$, every vertex in $Z_{i,s+1}$ has a non-neighbour in $X_{s+1}$. In case 1, 
this is true since $Z_{i,s+1}=\emptyset$, so we may assume that case 2 applies. In the notation of case 2,
we recall that $Z_{i,s+1}=Y_{i,s}\setminus W$, and so every vertex in $Z_{i,s+1}$ has a non-neighbour in $X_{s+1}$ as required.
This proves {\bf(C3)}.

For {\bf(C4)}, it suffices to show that for $i<j< s+1$, $Z_{i,j}$ is anticomplete to $X_{s+1}\cup N_{s+1}$. But this is true
since $Z_{i,j}$ is anticomplete to $C$ and $X_{s+1}\cup N_{s+1}\subseteq C$. This proves {\bf(C4)}.

For {\bf(C5)}, we must show that for all $i<j\le s+1$, either
\begin{itemize}
\item some vertex in $X_j$ is anticomplete to $Y_{i,s+1}$, and $Z_{i,j}=\emptyset$, or
\item $X_j$ is complete to $Y_{i,s+1}$, and $Z_{i,j}$ covers $N_j$.
\end{itemize}
If $j\le s$, then the claim holds since $Y_{i,s+1}\subseteq Y_{i,s}$ and either 
\begin{itemize}
\item some vertex in $X_j$ is anticomplete to $Y_{i,s}$, and $Z_{i,j}=\emptyset$, or
\item $X_j$ is complete to $Y_{i,s}$, and $Z_{i,j}$ covers $N_j$.
\end{itemize}
Consequently we may assume that $j=s+1$. If case 1 applies, let $x$ be as in case 1; then $x$ is anticomplete to 
$Y_{i,s+1}$ from the definition of $Y_{i,s+1}$, and $Z_{i,s+1}=\emptyset$, so the claim holds. If case 2 applies,
let $W$ be as in case 2; then $X_j$ is complete to $Y_{i,s}$ since $Y_{i,s+1}=W$, and $Z_{i,s+1}$ covers $N_{s+1}$
since this was shown just before the definition of ``case 2''. This proves {\bf(C5)}, and so completes the proof
that $G$ admits an $h$-cable of length $s+1$ with chromatic number more than $\sigma_{s+1}$. 

We have shown then that for $s=0\l t$, $G$ admits an $h$-cable of length $s$ with chromatic number more than $\sigma_s$. In particular,
$G$ admits an $h$-cable of length $t$ with chromatic number more than $\sigma_t=c$. This proves \ref{getcable}.~\bbox

\bigskip

By combining \ref{usecable} and \ref{getcable}, we deduce:
\begin{thm}\label{mainthm2}
Let $k,\kappa,\tau,\ell\ge 0$ and $h\ge 1$, and let $\phi:\mathbb{N}\rightarrow\mathbb{N}$ be nondecreasing.
Then there exists $c(\tau)\ge 0$ with the following property.
Let $G$ be a graph such that:
\begin{itemize}
\item $G$ has no hole of length at least $\ell$;
\item $\omega(G)\le k$;
\item $\chi(H)\le \kappa$ for every induced subgraph $H$ of $G$ with $\omega(H)<k$; 
\item $G$ is $(h,\phi)$-clique-controlled; and
\item $\chi(N^2(X))\le \tau$ for every $(h+1)$-clique $X$ of $G$.
\end{itemize}
Then $\chi(G)\le c(\tau)$.
\end{thm}
Of course $c(\tau)$ depends on all of $k,\kappa,\tau,\ell,h$ and $\phi$, but this notation is convenient.

\section{Proof of the main theorem}

We need the following. A somewhat stronger version was proved in~\cite{pentagonal},
but we give a proof here to make the paper self-contained.

\begin{thm}\label{longholelemma}
Let $\ell\ge 4$, $\kappa\ge 0$ and $\tau\ge 0$ be integers, and let $G$ be a graph with no hole of length at least $\ell$, such that
$\chi(N^1(v))\le \kappa$ and $\chi(N^2(v))\le \tau$ for every vertex $v$.
Then $\chi(G)\le 2(\ell-3)(\kappa+\tau)+1$.
\end{thm}
\Proof
Let $G_1$ be a component of $G$ with $\chi(G_1) = \chi(G)$, let $z_0\in V(G_1)$, and for $i\ge 0$ let $L_i$ be the set of
vertices of $G_1$ with distance $i$ from $z_0$. Choose $k$ such that $\chi(L_k)\ge \chi(G_1)/2$. If $k=0$ then the theorem holds,
so we may assume that $k\ge 1$.
Let $C_0$ be the vertex set of a component
of $G[L_k]$ with maximum chromatic number. Choose $v_0\in L_{k-1}$ with a neighbour in $C_0$. Let $t=\ell-3$, and suppose that
$\chi(C_0)>t\kappa+t\tau$. We claim that :
\\
\\
(1) {\em For all $i$ with $0\le i\le t$, there is
an induced path $v_0\d v_1\c v_i$ where $v_1\l v_i\in C_0$, and a subset $C_i$ of $C_0$ such that $G[C_i]$ is connected,
$\chi(C_i)> (t-i)\kappa+t\tau$, $v_i$ has a neighbour in $C_i$, and $v_0\l v_{i-1}$ have no neighbours in $C_i$.}
\\
\\
For this is true when $i = 0$; suppose it is true for some value of $i<t$, and we prove it is also true for $i+1$.
Let $N$ be the set of neighbours of $v_i$ in $C_i$. Thus
$$\chi(C_i\setminus N)\ge \chi(C_i)-\kappa> (t-i-1)\kappa+t\tau\ge 0,$$
and so $C_i\setminus N\ne \emptyset$; let $C_{i+1}$ be the vertex set of a component of $G[C_i\setminus N]$ with
maximum chromatic number. Thus
$\chi(C_{i+1})> (t-i-1)\kappa+t\tau$.
Choose $v_{i+1}\in N$ with a neighbour in $C_{i+1}$. This completes the inductive definition of $v_1\l v_i$ and $C_i$,
and so proves (1).

\bigskip
In particular, such a path $v_0\c v_t$ and subset $C_t$ exist. Since $\chi(C_t)>t\tau$, there
is a vertex $v\in C_t$ in none of the sets $N^2_G(v_i)\;(0\le i\le t-1)$, and therefore with distance at least three from
all of $v_0\l v_{t-1}$, since $t\ge 1$. Choose
$u\in L_{k-1}$ adjacent to $v$; then $u$ has distance at least two from all of $v_0\l v_{t-1}$. Let $P$ be an induced path
of $G[C_t\cup \{u,v_t\}]$ between $u,v_t$; thus $P$ has length at least one. Let $Q$ be an induced path of $G$ between
$u,v_0$ with all internal vertices in $L_0\cup\cdots\cup L_{k-2}$; then $Q$ has length at least two. The union of
$P,Q$ and $v_0\d v_1\c v_t$ is a hole of length at least $t+3=\ell$, which is impossible.

This proves that $\chi(C_0)\le t\kappa+t\tau$. Consequently $\chi(L_k)\le t(\kappa+\tau)$, and so $\chi(G)\le 2t(\kappa+\tau)$.
This proves \ref{longholelemma}.~\bbox

\bigskip

From \ref{longholelemma} we deduce:
\begin{thm}\label{control1}
Let $\ell\ge 4$, and let $k\ge 1$ and $\kappa\ge 0$ be such that
$\chi(H)\le \kappa$
for every graph $H$ with no hole of length at least $\ell$ and $\omega(H)<k$.
For $x\ge 0$ let $\phi_1(x) = 2(\ell-3)(\kappa+x)+1$.
Then every graph $G$ with no hole of length at least $\ell$
and with $\omega(G)\le k$ is $(1,\phi_1)$-clique-controlled.
\end{thm}
\Proof
Let $G$ be a graph with no hole of length at least $\ell$ and with $\omega(G)\le k$. Let $n\ge 0$,
and let $H$ be an induced subgraph of $G$ with $\chi(H)>\phi_1(n)$. Consequently $V(H)\ne \emptyset$; choose
$v\in V(H)$
with $\chi(N^{2}_H(v))$ maximum, $\chi(N^{2}_H(v))=\tau$ say.
Since $H$ has no hole of length at least $\ell$, and
$\chi(N_H(u))\le \kappa$ and $\chi(N^2_H(u))\le \tau$ for every vertex $u$ of $H$, \ref{longholelemma} implies that
$\chi(H)\le 2(\ell-3)(\kappa+\tau)+1$,
and so $\phi_1(n)<\chi(H)\le \phi_1(\chi(N^{2}_H(v)))$. Consequently $\chi(N^{2}_H(v))>n$. This proves \ref{control1}.~\bbox

\bigskip

We claim:
\begin{thm}\label{controlh}
Let $\ell\ge 4$, and let $k\ge 1$ and $\kappa\ge 0$ be such that
$\chi(H)\le \kappa$
for every graph $H$ with no hole of length at least $\ell$ and $\omega(H)<k$.
For all $h$ with $1\le h\le k$ there is a nondecreasing function $\phi_h:\mathbb{N}\rightarrow \mathbb{N}$
such that every graph $G$ with no hole of length at least $\ell$
and with $\omega(G)\le k$ is $(h,\phi_h)$-clique-controlled.
\end{thm}
\Proof We proceed by induction on $h$. In view of \ref{control1}, the result holds for $h=1$, so we may assume that $h<k$
and the result holds for $h$, and we will prove it holds for $h+1$. Since the result holds for $h$, $\phi_h$ exists as in the theorem.
By \ref{mainthm2}, for each $\tau\ge 0$, there exists $c(\tau)$ as in \ref{mainthm2} with $\phi$ replaced by $\phi_h$. 
For each $n\ge 0$, let $\phi_{h+1}(n)=\max_{0\le \tau\le n}c(\tau)$; we claim that $\phi_{h+1}$ satisfies the theorem.
For let $G$ be a graph with 
no hole of length at least $\ell$, and 
$\omega(G)\le k$. It follows that $G$ is $(h,\phi_h)$-clique-controlled. We must show that $G$ is $(h+1,\phi_{h+1})$-clique-controlled.
Thus, let $H$ be an induced subgraph of $G$, and let $\chi(H)>\phi_{h+1}(n)$ for some $n\ge 0$; we must show that 
there is an $(h+1)$-clique $X$ of $H$ such that $\chi(N^2_H(X))>n$.
Let $\tau$ be the maximum of $\chi(N^2_H(X))$ over all $(h+1)$-cliques $X$ of $H$, or $0$ if there is no such $X$. By \ref{mainthm2},
$\chi(H) \le c(\tau)$, and so $c(\tau)>\phi_{h+1}(n)$. It follows that $\tau>n$, and so 
there is an $(h+1)$-clique $X$ of $H$ such that $\chi(N^2_H(X))>n$. This proves \ref{controlh}.~\bbox

\bigskip

\noindent{\bf Proof of \ref{mainthm}.\ \ }
By induction on  $k$, we may assume that there exists $\kappa\ge 0$ such that $\chi(H)\le \kappa$
for every graph $H$ with no hole of length at least $\ell$ and $\omega(H)<k$.
Given $k,\ell$, let $\phi_k$ be as in \ref{controlh}, and let $c=\phi_k(0)$. We claim that 
$c$ satisfies \ref{mainthm}. For let $G$ be a graph with no hole of length at least $\ell$
and with $\omega(G)\le k$; then $G$ is $(k,\phi_{k})$-clique-controlled, by \ref{controlh} with $h=k$. If $\chi(G)>\phi_k(0)$ 
then there is a $k$-clique $X$ of $G$ such that $\chi(N^2_G(X))>0$, which is impossible since 
$\omega(G)\le k$ and so $\chi(N^2_G(X))=0$ for every $k$-clique $X$.
This proves that $\chi(G)\le \phi_k(0)=c$, and so proves \ref{mainthm}.~\bbox

\end{document}